\documentclass[12pt]{article}
\def\q{\hfill\rule{1ex}{1ex}}
\def\0{\emptyset}

\def\q{\hfill\rule{1ex}{1ex}}

\def\n{\noindent}

\usepackage{graphicx}
\usepackage{authblk}
\usepackage{amsfonts}
\usepackage{enumitem}
\usepackage{calc}
\usepackage{amsthm}
\usepackage{amssymb}
\usepackage{amsmath}
\usepackage{amsthm}
\usepackage{tikz}
\usepackage{pgfplots}

\begin{document}
\title{\bf Minimal colorings for properly colored subgraphs in complete graphs}

\author[1,2,4]{{\small\bf Chunqiu Fang}\thanks{ fcq15@mails.tsinghua.edu.cn, supported in part by CSC(No. 201806210164).}}
\author[2]{{\small\bf Ervin Gy\H ori}\thanks{ gyori.ervin@renyi.mta.hu, supported in part by the National Research, Development and Innovation Office NKFIH, grants K116769, K117879 and K126853.}}
\author[2,3,5]{{\small\bf Jimeng Xiao}\thanks{ xiaojimeng@mail.nwpu.edu.cn, supported in part by CSC(No. 201706290171).}}

\affil[1]{Department of Mathematical Sciences, Tsinghua University, Beijing 100084, China }
\affil[2]{Alfr\'ed R\'enyi Institute of Mathematics, Hungarian Academy of Sciences, Re\'altanoda u.13-15, 1053 Budapest, Hungary}
\affil[3]{Department of Applied Mathematics, Northwestern Polytechnical University, Xi'an, China}
\affil[4]{Yau Mathematical Sciences Center, Tsinghua University, Beijing 10084, China}
\affil[5]{Xi'an-Budapest Joint Research Center for Combinatorics, Northwestern Polytechnical University, Xi'an, China}

\date{}

\maketitle\baselineskip 16.7pt

\begin{abstract} Let $pr(K_{n}, G)$ be the maximum number of colors in an edge-coloring of $K_{n}$ with no properly colored copy of $G$. In this paper, we show that
$pr(K_{n}, G)-ex(n, \mathcal{G'})=o(n^{2}), $
where $\mathcal{G'}=\{G-M: M \text{ is a matching of }G\}$. Furthermore, we determine the
value of $pr(K_{n}, P_{l})$ for $l\ge 27$ and $n\ge 2l^{3}$  and the exact value of $pr(K_{n}, G)$, where $G$ is $C_{5}, C_{6}$ and $K_{4}^{-}$, respectively.  Also, we give an upper bound and a lower bound of $pr(K_{n}, K_{2,3})$.

\end{abstract}

~\\
Keywords: Properly colored subgraph, Tur\'an number, Anti-Ramsey number.
~\\

\n{\large\bf 1. Introduction}
~\\

We call a subgraph of an edge-coloring graph {\em rainbow}, if all of its edges have different colors. While a subgraph is called {\em properly colored} (also can be called {\em locally rainbow}), if any two adjacent edges receive different colors. The {\em anti-Ramsey number} of a graph $G$ in a complete graph $K_{n}$, denoted by $ar(K_{n}, G)$, is the maximum number of colors in an edge-coloring of $K_{n}$ with no rainbow copy of $G$. Namely, $ar(K_{n}, G)+1$ is the minimum number $k$ of colors  such that any $k$-edge-coloring of $K_{n}$ contains a rainbow copy of $G$.
In this paper, we let $pr(K_{n}, G)$ be the maximum number of colors in an edge-coloring of $K_{n}$ with no properly colored copy of $G$. Namely, $pr(K_{n}, G)+1$ is the minimum number $k$ of colors  such that any $k$-edge-coloring of $K_{n}$ contains a properly colored copy of $G$.

Given a family $\mathcal{F}$ of graphs, we call a graph $G$ an {\em $\mathcal{F}$-free} graph, if $G$ contains no graph in $\mathcal{F}$ as a subgraph. The {\em Tur\' an number }$ex(n, \mathcal{F})$ is the maximum number of edges in a graph $G$ on $n$ vertices which is $\mathcal{F}$-free. Such a graph $G$ is called an {\em extremal} graph, and the set of extremal graphs is denoted by $EX(n, \mathcal{F})$. The celebrated result of  Erd\H os-Stone-Simonovits Theorem \cite{erdos3, erdos1} states that for any $\mathcal{F}$ we have
\begin{align} \tag{*}
ex(n, \mathcal{F})=(\frac{p-1}{2p}+o(1))n^{2},
\end{align}
where $p=\Psi(\mathcal{F})=\min\{\chi(F): F\in\mathcal{F}\}-1,$ is the {\em subchromatic number}.

The anti-Ramsey number was introduced by Erd\H os, Simonovits and S\' os in \cite{erdos2}. There they showed that $ar(K_{n}, G)-ex(n, \mathcal{G})=o(n^{2})$, where $\mathcal{G}=\{G-e: e\in E(G)\}$ and  by (*), they showed that $ar(K_{n}, G)=(\frac{d-1}{2d}+o(1))n^{2}$, where $d=\Psi(\mathcal{G})$. This determined $ar(K_{n}, G)$ asymptotically when $\Psi(\mathcal{G})\ge 2.$ In case $\Psi(\mathcal{G})= 1$, the situation is more complex. Already the cases when $G$ is a tree or a cycle are nontrival. For a path $P_{k}$ on $k$ vertices, Simonovits and S\' os \cite{Simonovits} proved $ar(K_{n}, P_{2t+3+\epsilon})=tn-\binom{t+1}{2}+1+\epsilon$, for large $n$, where $\epsilon =0 $ or $1$. Jiang \cite{Jiang1} showed $ar(K_{n}, K_{1,p})=\lfloor\frac{n(p-2)}{2}\rfloor+\lfloor\frac{n}{n-p+2}\rfloor$ or possibly this value plus one if certain conditions hold. For  a general tree $T$ of $k$ edges, Jiang and West \cite{Jiang2} proved $\frac{n}{2}\lfloor\frac{k-2}{2}\rfloor+O(1)\le ar(K_{n}, T)\le ex(n, T)$ for $n\ge 2k$ and conjectured that $ar(K_{n}, T)\le \frac{k-2}{2}n+O(1).$ For cycles, Erd\H os, Simonovits and S\' os \cite{erdos2} conjectured that for every fixed $k\ge 3$, $ar(K_{n}, C_{k})=(\frac{k-2}{2}+\frac{1}{k-1})n+O(1)$, and proved that for $k=3$. Alon \cite{Alon} proved this conjecture for $k=4$ and gave some upper bounds for $k\ge 5$.  Finally, Montellano-Ballesteros and Neumann-Lara \cite{Montellano2} completely proved this conjecture.
For cliques, Erd\H os, Simonovits and S\' os \cite{erdos2} showed $ar(K_{n}, K_{p+1})=ex(n, K_{p})+1$ for $p\ge 3$ and sufficiently large $n$.  Montellano-Ballesteros and Neumann-Lara \cite{Montellano} and independently Schiermeyer \cite{Schiermeyer}  showed that $ar(K_{n}, K_{p+1})=ex(n, K_{p})+1$ holds for every $n\ge p\ge3$. For complete bipartite graphs $K_{s,t}, s\le t$, Axenovich and Jiang \cite{Axenovich} showed that  $ar(K_{n}, K_{2,t})=ex(n, K_{2,t-1})+O(n)$. Krop and York \cite{Krop} showed that  $ar(K_{n}, K_{s,t})=ex(n, K_{s,t-1})+O(n)$. Also, there are many other results about anti-Ramsey number. We mention the excellent survey by Fujita, Magnant, and Ozeki \cite{Fujita} for more conclusions on this topic.

The minimum number of colors guaranteeing the existence of properly colored subgraphs in an edge-colored complete graph was studied by Manoussakis, Spyratos, Tuza and Voigt in \cite{Manoussakis}.  For cliques, they showed that

 \noindent{\bf Theorem 1.} (\cite{Manoussakis})  For $t\ge 3$, let $b=\lfloor\frac{t-1}{2}\rfloor$, we have  $pr(K_{n}, K_{t})=(\frac{b-1}{2b}+o(1))n^2.$

For paths and cycles, they showed  \cite{Manoussakis} that $pr(K_{n}, P_{n})=\binom{n-3}{2}+1$ for large $n$ and $pr(K_{n}, C_{n})=\binom{n-1}{2}+1$. Also, they gave a conjecture on cycles as follows.

\noindent{\bf Conjecture 1.} (\cite{Manoussakis}) Let $n>l\ge 4$. Assume that $K_{n}$ is colored with at least $k$ colors, where 
\begin{align} \nonumber
k=
 \begin{cases}
 \frac{1}{2}l(l+1)+n-l+1, \, \text{if} \, n\le \frac{10l^{2}-6l-18}{6(l-3)},\\
 \frac{1}{3}ln-\frac{1}{18}l(l+3)+2, \, \text{if}\, n\ge \frac{10l^{2}-6l-18}{6(l-3)}.\\
 \end{cases}
 \end{align}
 Then, $K_{n}$ admits a properly colored cycle of length $l+1$.

In this  paper, we generalize Theorem 1 to arbitrary graph $G$ which shows that $pr(K_{n}, G)$ is related to the Tur\' an number like the anti-Ramsey number:

\noindent{\bf Theorem 2.} Let $G$ be a graph and $\mathcal{G'}=\{G-M: M \text{ is a matching of }G\}$, then $pr(K_{n}, G)\ge ex(n, \mathcal{G'})+1$  and $pr(K_{n}, G)=(\frac{d-1}{2d}+o(1))n^2,$ where $d=\Psi(\mathcal{G'})$.

We will prove Theorem 2 in Section 2 by the method used in the proof  of Theorem 1 in \cite{Manoussakis}.  Theorem 2  determines $pr(K_{n}, G)$ asymptotically when $\Psi(\mathcal{G'})\ge2$. As the anti-Ramsey number, the case $\Psi(\mathcal{G'})= 1$ is more complex.

In Section 3, we will determine $pr(K_{n}, P_{l})$ for $l\ge 27$ and large $n$  by proving the following theorem.

\noindent{\bf Theorem 3.} Let $P_{l}$ be a path with  $l\ge 27$ and $l\equiv r_{l}(\mod3)$, where $0\le r_{l}\le 2$. For $n\ge 2l^{3}$, we have $pr(K_{n}, P_{l})= (\lfloor\frac{l}{3}\rfloor-1)n-\binom{\lfloor\frac{l}{3}\rfloor}{2}+1+r_{l}.$

Although we just prove Theorem 3 for $l\ge 27$, we are sure that it is also true for  $l\le 26$.

For cycles, we slightly improved the lower bound of Conjecture 1 (See Proposition 4.1). Also, We modify Conjecture 1 as follows.

\noindent{\bf Conjecture 2.} Let $C_{k}$ be a cycle on $k$ vertices and $(k-1)\equiv r_{k-1}$ $(\mod 3),$ where $ 0\le r_{k-1}\le 2$, for $n\ge k,$
$$pr(K_{n}, C_{k})=\max\biggm\{\binom{k-1}{2}+n-k+1, \lfloor \frac{k-1}{3} \rfloor n-\binom{\lfloor \frac{k-1}{3} \rfloor+1}{2}+1+r_{k-1}\biggm\}.$$
 
 It is easy to see that $pr(K_{n}, C_{3})=ar(K_{n}, C_{3})=n-1$. Li, Broersma and Zhang \cite{Li}, and later Xu, Magnant, and Zhang\cite{Xu} showed that $pr(K_{n}, C_{4})=n$ for $n\ge4$. We consider $C_{5}$ and $C_{6}$ in section 4.

\noindent{\bf Theorem 4.} (a) $pr(K_{n}, C_{5})=n+2$ for $n\ge5$;

(b) $pr(K_{n}, C_{6})=n+5$ for $n\ge6$.

In Section 5, we consider two graphs $K_{4}^{-}$  and $K_{2,3}$, where $K_{4}^{-}$ is the graph $K_{4}$  minus one edge of it.

\noindent{\bf Theorem 5.} For $n\ge 3$, $pr(K_{n}, K_{4}^{-})= \lfloor\frac{3(n-1)}{2}\rfloor. $

\noindent{\bf Theorem 6.} For $n\ge 5$,  $\frac{7}{4}n+O(1)\le pr(K_{n}, K_{2,3})\le 2n-1.$

~\\

{\bf Notations:}  Let $G $ be a simple undirected graph. For $x\in V (G)$, we denote the {\em neighborhood} and the {\em degree} of $x$ in $G$ by $N_G(x)$ and $d_G(x)$, respectively.  The {\em maximum degree}  of $G$ is denoted by  $\Delta(G)$. The {\em common neighborhood} of $U\subset V(G)$ is the set of vertices in $V(G)\setminus U$ that are adjacent to every vertex in $U$.
We will use $G- x$ to denote the graph that arises from $G$ by deleting the vertex $x\in V (G)$. For $\emptyset\not= X \subset V (G)$, $G[X]$ is the subgraph of $G$ induced by $X$ and $G -X$ is the subgraph of $G$ induced by $V (G)\setminus X$. Given a graph $G= (V, E)$, for any (not necessarily disjoint) vertex sets $A, B\subset V$, we let $E_{G}(A, B):=\{uv\in E(G)| u\ne v, u\in A, v\in B\}$.  We  use $\overline{G}$ to denote the complement of $G$. Given two vertex disjoint graphs $G_{1}$ and $G_{2}$, we denote by $G_{1}+G_{2}$  the {\em join} of graphs $G_{1}$ and $G_{2}$, that is the graph obtained from $G_{1}\cup G_{2}$ by joining each vertex of $G_{1}$ with each vertex of $G_{2}$.

Given an edge-coloring $c$ of $G$, we denote the color of an edge $uv$ by $c(uv)$.  A color $a$ is  {\em starred} (at $x$) if all the edges with color $a$ induce a star $K_{1,r}$ (centered at the vertex $x$). We let $d^{c}(v)=|\{a\in C(v) | a$ is starred at $ v\}|$. For a subgraph $H$ of $G$, we denote $C(H)=\{c(uv)|\, uv\in E(H)\}$.
A {\em representing subgraph} in an edge-coloring of $K_n$ is a spanning subgraph
containing exactly one edge of each color.

\vskip.2cm
\n{\large\bf 2. Proof of Theorem 2}

\vskip.2cm
In this section, we will prove Theorem 2 by a similar argument used in the proof of Theorem 1 in \cite{Manoussakis}.

\noindent{\bf Theorem 2.} Let $G$ be a graph and $\mathcal{G'}=\{G-M: M \text{ is a matching of }G\}$, then $pr(K_{n}, G)\ge ex(n, \mathcal{G'})+1$  and $pr(K_{n}, G)=(\frac{d-1}{2d}+o(1))n^2,$ where $d=\Psi(\mathcal{G'})$.

\begin{proof}
 Let $F$ be a graph in $EX(n, \mathcal{G'})$. We color the edges of $K_{n}$ as follows. Take a subgraph $F$ of $K_{n}$, and assign distinct colors to all of $E(F)$ and a new color $c_{0}$ to all the remaining edges.  Suppose there is a properly colored $G$, then $M=\{e\in E(G), e \text{ is colored with } c_{0}\}$ is a matching of $G$, and $G-M\subset F$. By the definition of $\mathcal{G'}$, we have $G-M \in \mathcal{G'}$, and this is a contradiction with $F$ being $\mathcal{G'}$-free. Thus we have $pr(K_{n}, G)\ge ex(n, \mathcal{G'})+1=(\frac{d-1}{2d}+o(1))n^2$ by (*).

  Let $G_{0}=G-M_{p}$, where $M_{p}$ is a $p$-matching of $G$ and $\chi(G_{0})=d+1$. We prove that for every fixed $\varepsilon > 0$,  and for $n$ large enough with respect to $n_{0}=|V(G)|$ and $\varepsilon$, there is a properly colored copy of $G$ in any $(\frac{d-1}{2d}+\varepsilon)n^2$-edge coloring of $K_{n}$. In a representing subgraph of $K_{n}$ with $(\frac{d-1}{2d}+\varepsilon)n^2$ edges,  for an arbitrarily fixed $s$, and for $n$ sufficiently large, by (*), there exists a complete $(d+1)$-partite subgraph $K_{s,s,\cdots,s}$ with $s$ vertices in each class. We take $s=2^{n_{0}+d+1}.$

Denote by $V$ the vertex set of $K_{s, s, \cdots, s}$ and by $V_{1}, V_{2}, \cdots, V_{d+1}$ its vertex classes. We apply the procedure that follows.

For each $i=1, 2, \cdots, d+1$ do sequentially the following:
\vskip.2cm

\noindent (1) Select arbitrarily $2^{n_{0}+d+1-i}$ pairwise disjoint pairs $\{u_{ij}, v_{ij}\}$ in $V_{i}$, $ j= 1,2\cdots, 2^{n_{0}+d+1-i}.$

\noindent (2) For $j=1, 2, \cdots, 2^{n_{0}+d+1-i}$, delete from $K_{s, s, \cdots, s}$ the (at most one) vertex $z\in V\setminus V_{i}$ for which either $c(zu_{ij})=c(u_{ij}v_{ij})$ or $c(zv_{ij})=c(u_{ij}v_{ij})$, and if $z$ has already been selected in a previous pair $\{u_{i'j'}, v_{i'j'}\}$, for some $i'<i$, then also delete the other member of its pair.
\vskip.2cm

{\bf Claim 1.} It is possible to carry out the above procedure and that at the end of the execution, in each $V_{i}$, at least $2^{n_{0}}$ pairs remains undeleted.

{\bf The proof of Claim 1.} In the beginning, $V_{i}$ contains $2^{n_{0}+d+1}$ vertices, $i=1, 2,
\cdots, d+1$. In the first iteration, $i=1$, we can carry out (1) and (2) easily. Suppose we have carried out up to the $(i-1)$-st iteration. Before executing the $i$-th iteration observe that at most $\sum_{1\le j\le i-1}2^{n_{0}+d+1-j}= 2^{n_{0}+d+1}-2^{n_{0}+d+2-i}$ vertices have been deleted from $V_{i}$. Thus, $V_{i}$ contains at least $2^{n_{0}+d+2-i}$ vertices and it is enough to execute instruction (1) in the $i$-th iteration.

On the other hand, for any $i=1, 2, \cdots, d$, from the $(i+1)$-st iteration up to the end, due to instructions of type (2), at most  $\sum_{i+1\le j\le d+1}2^{n_{0}+d+1-j}=2^{n_{0}+d+1-i}-2^{n_{0}} $ pairs in $V_{i}$ have been delete and thus at least $2^{n_{0}}$ pairs in $V_{i}$ remains undeleted. Note also that $V_{d+1}$ contains $2^{n_{0}}$ pairs of vertices and there is no deletion of pair in $V_{d+1}$.\q

For $1\le i\le d+1$, let $\{x_{ij}y_{ij}: 1\le j\le 2^{n_{0}}\}$ be the $2^{n_{0}}$ pairs in $V_{i}$ which remains undeleted and $V'_{i}=\{x_{ij}, y_{ij}: 1\le j\le 2^{n_{0}}\}$. Let $$H=K_{s, s, \cdots, s}[V'_{1}\cup\cdots\cup V'_{d+1}]\bigcup (\cup_{i=1}^{d+1}\{x_{ij}y_{ij}: 1\le j\le 2^{n_{0}}\}). $$
 $H$ is properly colored, by Claim 1. Since $G_{0}= G-M_{p}$ and $\chi(G_{0})=d+1$, we have $H\supset G$. Thus $pr(K_{n}, G)\le (\frac{d-1}{2d}+o(1))n^2$.
\end{proof}


\vskip.2cm

\n{\large\bf 3. Paths}
\vskip.2cm

In this section, we consider the minimum number of colors guaranteeing the existence of properly edge-colored paths in an edge-colored complete graph, and prove Theorem 3. Before doing so, we will determine $pr(K_{n}, P_{l})$ for small $l$.

\noindent{\bf Proposition 3.1} (a) $pr(K_{n}, P_{3})=1,$ for $ n\ge 3$.

(b) $pr(K_{n}, P_{4})=2,$ for $ n\ge 4$.

(c) $pr(K_{n}, P_{5})=3,$ for $ n\ge 5$.

(d) $pr(K_{n}, P_{6})=n,$ for $ n\ge 6$.

\begin{proof}
(b)  Choose a vertex $v$ of $K_{n}$,  color all  edges incident to $v$ with color $c_{1}$ and color  all the remaining edges with color $c_{2}$. We use two colors and there is no properly colored $P_{4}$. Hence $pr(K_{n}, P_{4})\ge 2$.

 For $n\ge 5$, we have $pr(K_{n}, P_{4})\le ar(K_{n}, P_{4})=2$ (see \cite{Bialostocki}). Consider a $3$-edge-coloring of $K_{4}$. Let $V(K_{4})=\{u, v, x, y\}$. Then there is at least one edge in $E(\{u,v\}, \{x,y\})$, say $ux$ such that $c(ux)\ne c(uv)$ and $c(ux)\ne c(xy)$. Thus $vuxy$ is a properly colored $P_{4}$. Hence $pr(K_{n}, P_{4})\le 2$.

(c)  Choose  $u, v\in V(K_{n})$, assign $c_{1}$ to all edges incident with $u$, $c_{2}$ to all  edges incident with $v$ (except the edge $uv$) and  $c_{3}$ to all the remaining edges. We use three colors and there is no properly colored $P_{5}$. Hence $pr(K_{n}, P_{5})\ge 3$.

Consider a $4$-edge-coloring of $K_{n}, n\ge 5$, there is always a rainbow $P_{4}=u_{1}u_{2}u_{3}u_{4}$ since $ar(K_{n}, P_{4})=2$ when $ n\ge 5$. Since $|C(P_{4})|=|E(P_{4})|=3$, there is a color $c_{0}\notin C(P_{4})$.

Suppose there is no properly colored $P_{5}$ in the $4$-edge-coloring of $K_{n}$. Then for all $u\in V(K_{n})\backslash V(P_{4})$, it must be $c(uu_{1})=c(u_{1}u_{2})$, $c(uu_{4})=c(u_{3}u_{4})$, $c(uu_2)\in \{c(u_1u_2),c(u_2u_3)\}$ and $c(uu_3)\in \{c(u_2u_3),c(u_3u_4)\}$.
If $c(u_{1}u_{4})=c_{0},$ then $uu_{1}u_{4}u_{3}u_{2}$ is a properly colored $P_{5}$, a contradiction. If $c(u_{1}u_{3})=c_{0}$ or $c(u_{2}u_{4})=c_{0}$, say $c(u_{1}u_{3})=c_{0}$, then $u_{4}uu_{1}u_{3}u_{2}$ is a properly colored $P_{5}$, a contradiction. So there are two vertices $x,y\in V(K_{n})-V(P_{4})$ such that $c(xy)=c_{0}$. In this case, we have $u_{4}yxu_{2}u_{1}$ or $u_{4}yxu_{2}u_{3}$ is a properly colored $P_{5}$, a contradiction. Hence $pr(K_{n}, P_{5})\le 3$.








(d)  Choose a vertex $v$ of $K_{n}$, assign distinct colors to all the edges incident with vertex $v$ and a new color to all the remaining edges. We use $n$ colors and there is no properly colored $P_{6}$. Hence $pr(K_{n}, P_{6})\ge n$.

 Consider a $(n+1)$-edge-coloring of $K_{n}$. Then there is always a rainbow $P_{5}=u_{1}u_{2}u_{3}u_{4}u_{5}$ since $ar(K_{n}, P_{5})=n$ (see \cite{Bialostocki}).

Suppose there is no properly colored $P_{6}$ in the $(n+1)$-edge-coloring of $K_{n}$. Then
for all $u\in V(K_{n})\backslash V(P_{5})$, it must be $c(uu_{1})=c(u_{1}u_{2})$, $c(uu_{5})=c(u_{4}u_{5})$,  $c(uu_{2})\in \{c(u_{1}u_{2}), c(u_{2}u_{3})\}$ and $c(uu_{4})\in \{c(u_{3}u_{4}), $  $ c(u_{4}u_{5})\}$.
If there is a vertex $u\in V(K_{n})\backslash V(P_{5})$ such that $c(uu_{2})=c(u_{1}u_{2})$ and $c(uu_{4})=c(u_{4}u_{5})$, then at least one of $uu_{2}u_{3}u_{4}u_{5}u_{1}$ and $uu_{4}u_{3}u_{2}u_{1}u_{5}$ is a properly colored $P_{6}$ whatever $c(u_{1}u_{5})$ is, a contradiction. Hence  for all $u\in V(K_{n})\backslash V(P_{5})$,  if $c(uu_{2})=c(u_{1}u_{2})$ (resp. $c(uu_{4})=c(u_{4}u_{5})$), then  $c(uu_{4})\not=c(u_{4}u_{5})$ (resp. $c(uu_{2})\not=c(u_{1}u_{2})$).

If $c(u_{1}u_{5})\notin C(P_{5}),$ take $u\in V(K_{n})\backslash V(P_{5})$, then $uu_{1}u_{5}u_{4}u_{3}u_{2}$ is a properly colored $P_{6}$, a contradiction.

If $c(u_{1}u_{4})\notin C(P_{5})$ or $c(u_{2}u_{5})\notin C(P_{5})$,  say $c(u_{1}u_{4})\notin C(P_{5})$,
 take $u\in V(K_{n})\backslash V(P_{5})$, then $u_{5}uu_{1}u_{4}u_{3}u_{2}$ is a properly colored $P_{6}$, a contradiction.

Suppose $c(u_{1}u_{3})\notin C(P_{5})$ or $c(u_{3}u_{5})\notin C(P_{5})$,  say $c(u_{1}u_{3})\notin C(P_{5})$. Take $u\in V(K_{n})\backslash V(P_{5})$. If $c(uu_{2})=c(u_{2}u_{3})$, then $uu_{2}u_{1}u_{3}u_{4}u_{5}$ is a properly colored $P_{6}$, a contradiction. So $c(uu_{2})=c(u_{1}u_{2})$ which implies  $c(uu_{4})=c(u_{3}u_{4})$. Hence, $u_{1}u_{3}u_{2}uu_{4}u_{5}$ is a properly colored $P_{6}$, a contradiction.

If $c(u_{2}u_{4})\notin C(P_{5})$, take $u\in V(K_{n})\backslash V(P_{5})$, then at least one of $u_{1}uu_{3}u_{2}u_{4}u_{5}$ and $u_{1}u_{2}u_{4}u_{3}uu_{5}$ is a properly colored $P_{6}$ whatever $c(u_{1}u_{3})$ is, a contradiction.

Suppose there is   $u\in V(K_{n})\backslash V(P_{5})$ such that $c(uu_{3})\notin C(P_{5})$.
Then at least one of $u_{1}u_{2}uu_{3}u_{4}u_{5}$ and $u_{1}u_{2}u_{3}uu_{4}u_{5}$ is a properly colored $P_{6}$, a contradiction.

Since $n\ge 6$, there are two vertices $x,y\in V(K_{n})-V(P_{5})$ such that $c(xy)\notin C(P_{5})$.
Note that $c(xu_{3})\in C(P_{5})$. Then at least one of $u_{1}u_{2}u_{3}xyu_{5}$ and $u_{1}yxu_{3}u_{4}u_{5}$ is a properly colored $P_{6}$ whatever $c(xu_{3})$ is, a contradiction. Hence $pr(K_{n}, P_{6})\le n$.
\end{proof}

Here, we give the  lower bound of $pr(K_{n}, P_{l})$ by the following proposition.

\noindent{\bf Proposition 3.2} Let $n\ge l$ and $P_{l}$ be a path with  $l\equiv r_{l}$  $(\mod3), 0\le r_{l}\le 2$. We have
$$pr(K_{n}, P_{l})\ge\max\biggm\{\binom{l-3}{2}+1 , (\lfloor\frac{l}{3}\rfloor-1) n-\binom{\lfloor\frac{l}{3}\rfloor}{2}+1+ r_{l}\biggm\}.$$

\begin{proof}
We color the edges of $K_{n}$ as follows. For the first lower bound, we choose a $K_{l-3}$  and color it rainbow, and use one extra color for all the remaining edges. In such way, we use exactly $\binom{l-3}{2}+1$ colors and do not obtain a properly colored $P_{l}$.

 For the second lower bound,  we partition $K_{n}$ into two graphs $K_{\lfloor\frac{l}{3}\rfloor-1}+\overline{K}_{n-\lfloor\frac{l}{3}\rfloor+1}$ and $K_{n-\lfloor\frac{l}{3}\rfloor+1}$. First we color $K_{\lfloor\frac{l}{3}\rfloor-1}+\overline{K}_{n-\lfloor\frac{l}{3}\rfloor+1}$ rainbow and  color  $K_{n-\lfloor\frac{l}{3}\rfloor+1}$ with $(1+ r_{l})$ new colors without producing a properly colored $P_{3+ r_{l}}$ (See the proof of Proposition 3.1).
In such way, we use exactly $(\lfloor\frac{l}{3}\rfloor-1) n-\binom{\lfloor\frac{l}{3}\rfloor}{2}+1+ r_{l}$ colors and do not obtain a properly colored $P_{l}$.
\end{proof}

The proof of the following lemma is trivial. We will use it to prove Theorem 3.

\noindent{\bf Lemma 3.3} Let $P_{l}$ be a path with $l$ vertices, and $l\equiv r_{l}(\mod3), 0\le r_{l}\le 2$.  If an  edge-coloring of $K_{n}$ contains a rainbow copy of $K_{\lfloor\frac{l}{3}\rfloor-1, 2\lfloor\frac{l}{3}\rfloor+3}$ but does not contain a properly colored $P_{l}$, then it is the following coloring: denote by  $Q$ the vertices of $K_{n}-K_{\lfloor\frac{l}{3}\rfloor-1, 2\lfloor\frac{l}{3}\rfloor+3}$, by $X$ the smaller class of $K_{\lfloor\frac{l}{3}\rfloor-1, 2\lfloor\frac{l}{3}\rfloor+3}$ and by $Y$ the other one. Then $|C(K_{n}[Y])|\le 1+r_{l}$. Also, we have $|C(K_{n}[Y])\cup C(E_{K_{n}}(Y, Q))|\le 1+r_{l}$ and $|C(K_{n}[Y\cup Q])|\le 1+r_{l}$.  We get the most colors if the colors of all the edges between $X$ and $Y\cup Q$ and all the edges in $X$ are different, they differ from all the other edges and we use exactly $1+r_{l}$ colors in  $Y\cup Q$ such that there is no properly colored $P_{3+r_{l}}$ in $Y\cup Q$. Then the number of colors is $$(\lfloor\frac{l}{3}\rfloor-1)n-\binom{\lfloor\frac{l}{3}\rfloor}{2}+1+r_{l}.$$

Now, we will prove Theorem 3, and the idea comes from \cite{Simonovits}.

\noindent{\bf Theorem 3.} Let $P_{l}$ be a path with  $l\ge 27$. Let $l\equiv r_{l}(\mod3), 0\le r_{l}\le 2$. For $n\ge 2l^{3}$,  we have $pr(K_{n}, P_{l})= (\lfloor\frac{l}{3}\rfloor-1)n-\binom{\lfloor\frac{l}{3}\rfloor}{2}+1+r_{l}.$

\begin{proof}
We just need prove the upper bound. We shall use the following results of Erd\H os and Gallai (see \cite{erdos0}):

$$\text{(a)} \,~~~~~~~~~~~~~~~~~~~~~~~~~ex(n, P_{r})\le \frac{r-2}{2}n;$$
$$\text{(b)}\,~~~~~~~ ex(n,\{C_{r+1}, C_{r+2},\cdots\} )\le \frac{r(n-1)}{2}.$$

Consider an edge-coloring of $K_{n}$ using $pr(K_{n}, P_{l})$ colors without producing a properly colored $P_{l}$. Take a longest properly colored path $P_{s}=v_{1}v_{2}\cdots v_{s}$, where $s\le l-1.$ Denote by $G$ the graph obtained by choosing one edge from each remaining color such that the number of edges joining $P_{s}$ to the remaining $n-s$ vertices as large as possible. We would partition $V(G)\backslash V(P_{s})$ into three sets $U_{1}, U_{2}$ and $U_{3}$ as follows:

$U_{1}$ is the vertex set of  $V(K_{n})\backslash V(P_{s})$ not jointed to $P_{s}$ at all: neither by edges nor by paths;

$U_{2}$ is the set of isolated vertices of  $V(K_{n})\backslash V(P_{s})$ jointed to $P_{s}$ by edges;

$U_{3}=V(K_{n})\backslash (V(P_{s})\cup U_{1}\cup U_{2}).$

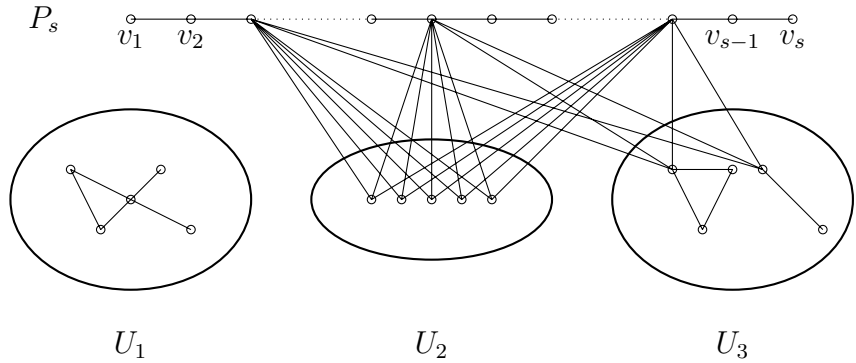
\begin{center}

\begin{tikzpicture}[scale=.8]

\draw[,thick] (0,0) arc(360:0:2cm and 1.5cm)  (5,0) arc(360:0:2cm and 1cm) (10,0) arc(360:0:2cm and 1.5cm)(-2,-2)node[below]{$U_{1}$}(3,-2)node[below]{$U_{2}$}(8,-2)node[below]{$U_{3}$};

\draw (-3,0.5)--(-2,0)--(-1.5,0.5) circle (2pt) (-3,0.5) circle (2pt)-- (-2.5,-0.5) circle (2pt)-- (-2,0) circle (2pt) --(-1,-0.5) circle (2pt)
(2,0) circle (2pt) (2.5,0) circle (2pt) (3,0) circle (2pt) (3.5,0) circle (2pt) (4,0) circle (2pt)
(7,0.5) circle (2pt)--(7.5,-0.5) circle (2pt)--(8,0.5) circle (2pt)--(7,0.5) (8.5,0.5) circle (2pt)--(9.5,-0.5) circle (2pt);

\draw (-3,3) node[left]{$P_{s}$} (-2,3)node[below]{$v_{1}$} circle (2pt)--(-1,3)node[below]{$v_{2}$}  circle (2pt)--(0,3) circle (2pt) (2,3) circle (2pt)--(3,3) circle (2pt)--(4,3) circle (2pt)--(5,3) circle (2pt) (7,3) circle (2pt)--(8,3)node[below]{$v_{s-1}$}  circle (2pt)--(9,3)node[below]{$v_{s}$}  circle (2pt);

\draw [rounded corners,dotted] (0,3) -- (2,3) (5,3) -- (7,3);

\draw (2,0)--(0,3) (2,0)--(3,3) (2,0)--(7,3)
(2.5,0)--(0,3) (2.5,0)--(3,3) (2.5,0)--(7,3)
(3,0)--(0,3) (3,0)--(3,3) (3,0)--(7,3)
(3.5,0)--(0,3) (3.5,0)--(3,3) (3.5,0)--(7,3)
(4,0)--(0,3) (4,0)--(3,3) (4,0)--(7,3)
(7,0.5)--(0,3) (7,0.5)--(3,3) (7,0.5)--(7,3)
(8.5,0.5)--(0,3) (8.5,0.5)--(3,3) (8.5,0.5)--(7,3);

\end{tikzpicture}

\end{center}

{\bf Claim 1.}  $G[U_{1}]$ contains no $P_{\lceil \frac{s+1}{2}\rceil}$, and $|E(G[U_{1}])|\le\frac{l-3}{4}|U_{1}|.$

{\bf Proof of Claim 1} Suppose $P_{\lceil \frac{s+1}{2}\rceil}=u_{1}u_{2}\cdots u_{\lceil \frac{s+1}{2}\rceil}$ is a path in $G[U_{1}]$. By the constructing of $G$, $c(u_{\lceil \frac{s+1}{2}\rceil}v_{\lceil\frac{s}{2}\rceil})\notin C(P_{\lceil \frac{s+1}{2}\rceil})$. Since $c(v_{\lceil\frac{s}{2}\rceil-1}v_{\lceil\frac{s}{2}\rceil})\ne c(v_{\lceil\frac{s}{2}\rceil}v_{\lceil\frac{s}{2}\rceil+1})$, at most one of $c(v_{\lceil\frac{s}{2}\rceil-1}v_{\lceil\frac{s}{2}\rceil})$ and $c(v_{\lceil\frac{s}{2}\rceil}v_{\lceil\frac{s}{2}\rceil+1})$ is the same as $c(u_{\lceil \frac{s+1}{2}\rceil}v_{\lceil\frac{s}{2}\rceil})$. So at least one of $u_{1}u_{2}\cdots u_{\lceil \frac{s+1}{2}\rceil}v_{\lceil\frac{s}{2}\rceil}\cdots v_{1}$ and $u_{1}u_{2}\cdots u_{\lceil \frac{s+1}{2}\rceil}v_{\lceil\frac{s}{2}\rceil}\cdots v_{s}$ is a properly colored path, a contradiction to the maximality of $P_{s}$. Hence, $G[U_{1}]$ contains no $P_{\lceil \frac{s+1}{2}\rceil}$. By (a), we have
$$|E(G[U_{1}])|\le \frac{1}{2}(\big\lceil \frac{s+1}{2}\big\rceil-2)|U_{1}|\le (\frac{1}{2}\big\lceil \frac{l}{2}\big\rceil-1)|U_{1}|\le\frac{l-3}{4}|U_{1}|.$$\q

{\bf Claim 2.} $E_{G}(U_{2}\cup U_{3}, \{v_{1}, v_{2}, v_{s-1}, v_{s}\})=\emptyset.$

{\bf Proof of Claim 2}
It is obvious that $E_{G}(U_{2}\cup U_{3}, \{v_{1}, v_{s}\})=\emptyset$ by the maximality of $P_{s}$.
Suppose that there is a vertex $u\in U_{2}\cup U_{3}$ such that $uv_{2}\in E(G)$ or $uv_{s-1}\in E(G)$, we say $uv_{2}\in E(G)$, then at least one of $uv_{1}v_{2}\cdots v_{s}$ and $v_{1}uv_{2}\cdots v_{s}$ is a properly colored path of order $s+1$, a contradiction to the maximality of $P_{s}$. \q

{\bf Claim 3.} $|E_{G}(U_{2}, P_{s})|\le (\lfloor\frac{l}{3}\rfloor-1)|U_{2}|.$

{\bf Proof of Claim 3}
For $v\in U_{2}$ and every three consecutive vertices $\{v_{i}, v_{i+1}, v_{i+2}\}\subset V(P_{s})$, we claim that $|E_{G}(v, \{v_{i}, v_{i+1}, v_{i+2}\})|\le1.$ Otherwise, at least two of $vv_{i}, vv_{i+1}, vv_{i+2}$ are edges of $G$. Then whatever $c(vv_{i})$ is, at least one of $v_{1}\cdots v_{i}vv_{i+1}v_{i+2}\cdots v_{s}$ and $v_{1}\cdots v_{i}v_{i+1}vv_{i+2}\cdots v_{s}$ is a properly colored path of order $s+1$, a contradiction to the maximality of $P_{s}$. By Claim 2, we have $|E_{G}(v, P_{s})|\le\lceil\frac{s-4}{3}\rceil\le\lceil\frac{l-5}{3}\rceil=\lfloor\frac{l}{3}\rfloor-1$ and $|E_{G}(U_{2}, P_{s})|\le (\lfloor\frac{l}{3}\rfloor-1)|U_{2}|.$\q

{\bf Claim 4.} $|E_{G}(U_{3}, P_{s})|+|E(G[U_{3}])|\le\frac{l+2}{4}|U_{3}|.$

{\bf Proof of Claim 4}
Take a component $H$ of $G[U_{3}]$ and let $r$ be the length of its longest cycle. If $H$ contains no cycles, then write $r=2$. For each vertex $u\in V(H)$, we can find a path $P_{r}\subset H$ starting from it. Hence, $E_{G}(u, \{v_{1},\cdots, v_{r},v_{s-r+1},\cdots, v_{s}\})=\emptyset$. Otherwise, we can find a properly colored path of order at least $s+1$. For any four consecutive vertices $\{v_{i}, v_{i+1}, v_{i+2}, v_{i+3}\}$, there are no two independent edges between $\{v_{i}, v_{i+1}, v_{i+2}, v_{i+3}\}$ and $V(H)$ in $G$. Otherwise, suppose $xv_{i}$ and $yv_{j}, j\in \{i+1, i+2, i+3\}$ are two independent edges between $\{v_{i}, v_{i+1}, v_{i+2}, v_{i+3}\}$ and $V(H)$ in $G$, where $x, y\in V(H)$. Take a path $P_{xy}$ of $H$ which connect $x$ and $y$. Then whatever $c(xv_{i+1})$ is, at least one of $v_{1}\cdots v_{i}xv_{i+1} \cdots v_{s}$ and $v_{1}\cdots v_{i}v_{i+1}xP_{xy}yv_{j}\cdots v_{s}$ is a prorperly colored path of order at least $s+1$,  a contradiction to the maximality of $P_{s}$. Hence, by (b), we have
$$|E_{G}(V(H), P_{s})|+|E(H)|\le \lceil\frac{s-2r}{4}\rceil|V(H)|+\frac{r|V(H)|-r}{2}\le \frac{s+3}{4}|V(H)|\le \frac{l+2}{4}|V(H)|.$$
By adding this up, we get $$|E_{G}(U_{3}, P_{s})|+|E(G[U_{3}])|\le \frac{l+2}{4}|U_{3}|.$$
\q

By Claims 1, 3 and 4, the number of colors is
\[\begin{split}pr(K_{n}, P_{l})&=|C(K_{n})| \le |C(P_{s})|+ |E(G)|\\
&\le \binom{s}{2}+|E(G[U_{1}])|+|E_{G}(U_{2}, P_{s})|+|E_{G}(U_{3}, P_{s})|+|E(G[U_{3}])|\\
&\le \binom{s}{2}+\frac{l-3}{4}|U_{1}|+(\lfloor\frac{l}{3}\rfloor-1)|U_{2}|+\frac{l+2}{4}|U_{3}|.\end{split}\]
Since $l\ge 27$, we have $\frac{l+2}{4}\le \lfloor\frac{l}{3}\rfloor-1-\frac{1}{4}.$ Let $U^{*}=\{u\in U_{2}: d_{G}(u)=\lfloor\frac{l}{3}\rfloor-1\}.$ Then
$$(\lfloor\frac{l}{3}\rfloor-1)n-\binom{\lfloor\frac{l}{3}\rfloor}{2}+1+r_{l}\le pr(K_{n},P_{l})\le \binom{s}{2}+(\lfloor\frac{l}{3}\rfloor-1-\frac{1}{4})(n-s-|U^{*}|)+(\lfloor\frac{l}{3}\rfloor-1)|U^{*}|.$$

Hence for $n\ge 2l^{3}$, we have $|U^{*}|\ge l^{3}$ and we can get at least $2\lfloor\frac{l}{3}\rfloor+3$ vertices $u_{1}, u_{2}, \cdots, u_{2\lfloor\frac{l}{3}\rfloor+3} \in U^{*}$ which have a common neighborhood of size $\lfloor\frac{l}{3}\rfloor-1$ in $G$. By Lemma 3.3, the proof is completed.
\end{proof}

\vskip.2cm

\n{\large\bf 4. Cycles}

\vskip.2cm

The lower bound of $pr(K_{n}, C_{k})$ was given roughly by Manoussakis, Spyratos, Tuza and Voigt in \cite{Manoussakis}. Here we prove the lower bound precisely again.

\noindent{\bf Proposition 4.1}  Let $C_{k}$ be a cycle on $k$ vertices and $k-1\equiv r_{k-1}$ $(\mod 3),$ where $ 0\le r_{k-1}\le 2$. For $n\ge k,$ we have
$$pr(K_{n}, C_{k})\ge\max\biggm\{\binom{k-1}{2}+n-k+1, \lfloor \frac{k-1}{3} \rfloor n-\binom{\lfloor\frac{k-1}{3}\rfloor+1}{2}+1+r_{k-1}\biggm\}.$$

\begin{proof}
We color the edges of $K_{n}$ as follows. For the first lower bound, we choose a $K_{k-1}$  and color it rainbow, and use one extra color for all the remaining edges. In such way, we use exactly $\binom{k-1}{2}+1$ colors and do not obtain a properly colored $C_{k}$.

 For the second lower bound,  we partition $K_{n}$ into two graphs $K_{\lfloor\frac{k-1}{3}\rfloor}+\overline{K}_{n-\lfloor\frac{k-1}{3}\rfloor}$ and $K_{n-\lfloor\frac{k-1}{3}\rfloor}$. First we color $K_{\lfloor\frac{k-1}{3}\rfloor}+\overline{K}_{n-\lfloor\frac{k-1}{3}\rfloor}$ rainbow and  color  $K_{n-\lfloor\frac{k-1}{3}\rfloor}$ with $(1+ r_{k-1})$ new colors without producing a properly colored $P_{3+ r_{k-1}}$ (See the proof of Proposition 3.1).
In such way, we use exactly $(\lfloor\frac{k-1}{3}\rfloor) n-\binom{\lfloor\frac{k-1}{3}\rfloor+1}{2}+1+ r_{k-1}$ colors and do not obtain a properly colored $C_{k-1}$.
\end{proof}

\noindent{\bf Conjecture 2.}  Let $C_{k}$ be a cycle on $k$ vertices and $(k-1)\equiv r_{k-1}$ $(\mod 3),$ where $ 0\le r_{k-1}\le 2$, for $n\ge k,$
$$pr(K_{n}, C_{k})=\max\biggm\{\binom{k-1}{2}+n-k+1, \lfloor \frac{k-1}{3} \rfloor n-\binom{\lfloor \frac{k-1}{3} \rfloor+1}{2}+1+r_{k-1}\biggm\}.$$

Although Li et al. \cite{Li} and later Xu et al.\cite{Xu} have got $pr(K_{n}, C_{4})=n$ for $n\ge 4$, here we will use a different method  to prove it again.
We denote a cycle $C_{k}$ with a pendant edge by $C_{k}^{+}$. Gorgol \cite{Gorgol} showed that $ar(K_{n}, C_{k}^{+})=ar(K_{n}, C_{k})$ for $n\ge k+1\ge 4$.  We will use this result  to prove the following proposition.

\noindent{\bf Proposition 4.2} (\cite{Li, Xu}) For $n\ge 4$, $pr(K_{n}, C_{4})=n$.

\begin{proof} By  Proposition 4.1, we have $pr(K_{n}, C_{4})\ge n$ for $n\ge 4.$ We will prove $pr(K_{n}, C_{4})\le n$  by induction  on $n$. The base case $n=4$ is obvious. For $n\ge 5$,  consider an $(n+1)$-coloring $c$ of $K_{n}$.  If there is a vertex $v$ such that $d^{c}(v)\le 1$, then $|C(K_{n}-v)|\ge n+1-1=(n-1)+1$ and there is a properly colored $C_{4}$ in $K_{n}-v$ by induction. Thus we assume that $d^{c}(v)\ge 2$, for all $v\in V$.
Since $ar(K_{n}, C_{3}^{+})=ar(K_{n}, C_{3})=n-1$, there is a rainbow $C_{3}^{+}$. Let the triangle be $xyzx$ and the pendant edge be $xu$. Let the edges $xy, yz, xz, xu$ have colors $1, 2, 3, 4$ respectively. We may assume $c(zu)\in \{ 2,4\}$; otherwise $xyzu$ is a properly colored $C_{4}$.

\begin{center}
\begin{tikzpicture}[scale=.8]

\draw (0,0)circle(2pt)node[left]{$y$} (0,2)circle(2pt)node[left]{$x$} (2,0)circle(2pt)node[below]{$z$} (2,2)circle(2pt)node[above]{$u$} (4,0)circle(2pt)node[below]{$v$}

(10,0)circle(2pt)node[left]{$y$} (10,2)circle(2pt)node[left]{$x$} (12,0)circle(2pt)node[below]{$z$} (12,2)circle(2pt)node[above]{$u$} (14,2)circle(2pt)node[below]{$v$} ;

\draw (0,0)--(0,2)--(2,2)--(2,0)--(0,2)  (0,0)--(2,0)--(4,0) (0,0)--(2,2)

(10,0)--(10,2)--(12,2)--(12,0)--(10,2)  (10,0)--(12,0) (12,2)--(14,2) (10,0)--(12,2);

\draw [rounded corners,dotted] (0,2) -- (4,0) (10,0) -- (14,2);

\draw (0,1)node{$1$} (1,0)node{$2$} (0.5,1.5)node{$3$}  (1,2)node{$4$} (0.5,0.5)node{$2$} (2,0.5)node{$2$} (3,0)node{$5$}

(10,1)node{$1$} (11,0)node{$2$} (10.5,1.5)node{$3$}  (11,2)node{$4$} (11.5,1.5)node{$4$} (12,0.5)node{$4$} (13,2)node{$5$};

\draw (2,-1)node{Case 1.} (12,-1)node{Case 2.};

\end{tikzpicture}
\end{center}

{\bf Case 1.} $c(zu)=2.$

We may assume that $c(yu)=2$; otherwise at least one of $xyuzx$ and $ xuyzx$ is a properly colored $C_{4}.$ Since $d^{c}(z)\ge 2$, there is a vertex $v\in V(K_{n})\setminus\{x, y, z, u\}$ such that $c(zv)$ is starred at $z$ and $c(zv)\ne 3$. Let $c(zv)=5$. Note that $c(xv)\ne 5$. Then at least one of $xyzvx$ and $xuzvx$ is a properly colored $C_{4}$.

{\bf Case 2.} $c(zu)=4.$

We may assume that $c(yu)=4$; otherwise at least one of $xyuzx$ and $ xuyzx$ is a properly colored $C_{4}.$ Since $d^{c}(u)\ge 2$, there is a vertex $v\in V(K_{n})\setminus\{x, y, z, u\}$ such that $c(uv)$ is starred at $u$ and $c(uv)\ne 4$. Let $c(zv)=5$. Note that $c(yv)\ne 5$. Thus at least one of $xuvyx$ and $zuvyz$ is a properly colored $C_{4}$.
\end{proof}

Now, we will use the same idea to prove Conjecture 1 for $k=5$.
Let $B$ be the {\em bull graph}, the unique graph on $5$ vertices with degree sequence $(1, 1, 2, 3, 3).$ Schiermeyer and Sot\' ak \cite{Schiermeyer} showed that $ar(K_{5}, B)=5$ and $ar(K_{n}, B)=n+1$ for $n\ge 6$. We will use this result  to prove the following proposition.

\noindent{\bf Proposition 4.3}  For $n\ge 5$, $pr(K_{n}, C_{5})=n+2$.

\begin{proof}
By  Proposition 4.1, we have $pr(K_{n}, C_{5})\ge n+2$ for $n\ge 5.$ We will prove $pr(K_{n}, C_{5})\le n+2$  by induction  on $n$. The base case $n=5$ is easy  since $pr(K_{n}, C_{n})=\binom{n-1}{2}+1$.  For $n\ge 6$, consider an $(n+3)$-edge-coloring $c$ of $K_{n}$.  If there is a vertex $v$ such that $d^{c}(v)\le 1$, then $|C(K_{n}-v)|\ge n+3-1=(n-1)+3$ and there is a properly colored $C_{4}$ by induction. Thus we assume that $d^{c}(v)\ge 2$, for all $v\in V$. Since $ar(K_{n}, B)=n+1$ for $n\ge 6$, there is a rainbow $B$. Let $E(B)=\{xy, xz, yz, yu, zv\}$ and the edges $xy, xz, yz, yu, zv$ have colors $1, 2, 3, 4, 5$ respectively. We can assume that $c(uv)\in \{4, 5\}$; otherwise $xyuvzx$ is a properly colored $C_{5}$. Assume, without loss of generality, that $c(uv)=4$. We may assume that $c(xu)\in \{1, 4\}$; otherwise $xyzvux$ is a properly colored $C_{5}$.

\begin{center}
\begin{tikzpicture}[scale=.8]

\draw (0,0)circle(2pt)node[left]{$z$}  (0,2)circle(2pt)node[left]{$y$}  (2,1)circle(2pt)node[right]{$x$}  (4,3)circle(2pt)node[above]{$u$}  (4,-1)circle(2pt)node[below]{$v$}  (6,2)circle(2pt)node[right]{$w$}

(10,0)circle(2pt)node[left]{$z$}  (10,2)circle(2pt)node[left]{$y$}  (12,1)circle(2pt)node[right]{$x$}  (14,3)circle(2pt)node[above]{$u$}  (14,-1)circle(2pt)node[below]{$v$}  ;

\draw (0,0)--(0,2)--(2,1)--(0,0)--(4,-1)--(4,3)--(2,1)  (4,3)--(0,2)--(4,-1) (4,3)--(6,2)--(0,2)

(10,0)--(10,2)--(12,1)--(10,0)--(14,-1)--(14,3)--(12,1)  (14,3)--(10,2);

\draw (0,1)node{$3$}  (1,1.5)node{$1$}  (1,0.5)node{$2$} (2,2.5)node{$4$} (2,-0.5)node{$5$} (4,1)node{$4$} (2.5,1.5)node{$1$}  (2,0.5)node{$4$} (5,2.5)node{$6$} (5,2)node{$3$}

(10,1)node{$3$}  (11,1.5)node{$1$}  (11,0.5)node{$2$} (12,2.5)node{$4$} (12,-0.5)node{$5$} (14,1)node{$4$} (13,2)node{$4$};

\draw [rounded corners,dotted] (0,0) -- (4,3)  (10,0) -- (14,3);

\draw (2,-2)node{Case 1.} (12,-2)node{Case 2.};

\end{tikzpicture}
\end{center}

{\bf Case 1.} $c(xu)=1$.

We may assume that $c(yv)=4$; otherwise at least one of $yvzxuy$ and $yvuxzy$ is a properly colored $C_{5}.$  Since $d^{c}(u)\ge 2$,  there is a vertex $w\in V(K_{n})\setminus \{x, y, z, u, v\}$ such that $c(uw)$ is starred at $u$ and $c(uw)\ne c(uz)$. Let $c(uw)=6$. Note that $c(yw)\ne 6$. We may assume that $c(yw)=3$;  otherwise $yzxuwy$  is a properly colored $C_{5}.$ Since $c(zu)\ne 6$, at least one of $uzvywu$ and $uzxywu$ is a properly colored $C_{5}$.

{\bf Case 2.} $c(xu)=4.$

Since $d^{c}(u)\ge 2,$ we consider the following two subcases.


\begin{center}
\begin{tikzpicture}[scale=.8]

\draw (0,0)circle(2pt)node[left]{$z$}  (0,2)circle(2pt)node[left]{$y$}  (2,1)circle(2pt)node[right]{$x$}  (4,3)circle(2pt)node[above]{$u$}  (4,-1)circle(2pt)node[below]{$v$}  (6,2)circle(2pt)node[right]{$w$}

(10,0)circle(2pt)node[left]{$z$}  (10,2)circle(2pt)node[left]{$y$}  (12,1)circle(2pt)node[right]{$x$}  (14,3)circle(2pt)node[above]{$u$}  (14,-1)circle(2pt)node[below]{$v$} (16,0)circle(2pt)node[right]{$w$} ;

\draw (0,0)--(0,2)--(2,1)--(0,0)--(4,-1)--(4,3)--(2,1)   (0,2)--(4,3)--(6,2)

(10,0)--(10,2)--(12,1)--(10,0)--(14,-1)--(14,3)--(12,1)  (14,3)--(10,2)
(10,0)--(14,3)  (10,2)--(14,-1)--(12,1) (14,-1)--(16,0);

\draw (0,1)node{$3$}  (1,1.5)node{$1$}  (1,0.5)node{$2$} (2,2.5)node{$4$} (2,-0.5)node{$5$} (4,1)node{$4$} (3,2)node{$4$}   (5,2.5)node{$$}

(10,1)node{$3$}  (11,1.5)node{$1$}  (11,0.5)node{$2$} (12,2.5)node{$4$} (12,-0.5)node{$5$} (14,1)node{$4$} (13,2)node{$4$}
(12,0.5)node{$1$} (12,1.5)node{$6$} (13,0)node{$1$} (15,-0.5)node{$7$};

\draw [rounded corners,dotted] (0,0) -- (6,2) (14,3)--(16,0);

\draw (2,-2)node{Case 2.1} (12,-2)node{Case 2.2};

\end{tikzpicture}
\end{center}


{\bf Case 2.1} There is a vertex $w\in V(K_{n})\setminus \{x, y, z, u, v\}$ such that $c(uw)$ is starred at $u$ and $c(uw)\ne 4$.

In this case, $c(uw)\ne c(zw)$. Then at least one of $zwuxyz$ and $zwuyxz$ is a properly colored $C_{5}$.

{\bf Case 2.2} $c(uz)$ is starred at $u$ and $c(uz)\ne 4$.

Let $c(uz)=6$. We can assume that $c(xv)=1$; otherwise at least one of $xvuzyx $ and $xvzuyx $  is a properly colored $C_{5}$. Also, we can assume that $c(yv)=1$; otherwise at least one of $yvuzxy$ and $yvzuxy$  is a properly colored $C_{5}$. Since $d^{c}(v)\ge 2$, there is a vertex $w\in V(K_{n})\setminus \{x, y, z, u, v\}$ such that $c(vw)$ is starred at $v$ and $c(vw)\ne 5$. We assume that $c(vw)=7$. Note that $c(uw)\ne 7$. Then at least one of $wuyzvw$ and $wuzxvw$ is a properly colored $C_{5}$.
\end{proof}


For $C_{6}$, we consider more cases to prove it.

\noindent{\bf Proposition 4.4} For $n\ge 6$, $pr(K_{n}, C_{6})=n+5$.

\begin{proof}
By Proposition 4.1, we have $pr(K_{n}, C_{6})\ge n+5$ for $n\ge 6.$ We will prove $pr(K_{n}, C_{6})\le n+5$ by induction on $n$. For $n=6$, $pr(K_{6}, C_{6})=\binom{6-1}{2}+1=11$. For $n=7$, $pr(K_{7}, C_{6})\le ar(K_{7}, C_{6})=12$ (see \cite{Montellano2}).
For $n\ge 8$, consider an $(n+6)$-edge-coloring $c$ of $K_{n}$. If there is a vertex $v$ such that $d^{c}(v)\le 1$, then $|C(K_{n}-v)|\ge n+6-1=(n-1)+6$ and there is a properly colored $C_{6}$ by induction. Thus we assume that $d^{c}(v)\ge 2$ for all $v\in V(K_{n})$. Let $G$ be a subgraph of $K_{n}$ such that $e\in E(G)$ if and only if the color $c(e)$ appears only once in $K_{n}$. We have $|E(G)|\ge 2n-(n+6)=n-6\ge 2$.

{\bf Case 1.} $\Delta(G)\ge 2 $.

In this case, $G$ contains a path of order 3. Let $P_{3}=xyz$ and $U=V(K_{n})\setminus \{x, y, z\}.$ For all $v\in U$ and all starred color $c_{v}$ at $v$, we take an edge with color $c_{v}$ to obtain a subgraph $H$ of $K_{n}$. Choose $H$  such that $|E_{H}(\{x, y, z\}, U)|$ as large as possible. 

{\bf Case 1.1}  $|E(H[U])|\ge 2$.

Let $u_{1}u_{2}, v_{1}v_{2}\in E(H[U])$. If $u_1\in\{v_1,v_2\}$ or $u_2\in\{v_1,v_2\}$, say $u_2=v_1$, then $c(xu_1)\ne c(u_1v_1)$ and $c(zv_2)\ne c(v_1v_2)$ by the choice  of $H$. Thus $xyzv_2v_1u_1x$ is a properly colored $C_{6}$. Now suppose $u_1u_2$ and $v_1v_2$ are two independent edges of $H$.  Assume that $c(u_{1}u_{2})$ and $c(v_{1}v_{2})$ are starred at $u_{1}, v_{1}$ respectively. Thus $c(u_{2}v_{2})\ne c(u_{1}u_{2})$  and $c(u_{2}v_{2})\ne c(v_{1}v_{2})$. By the choice  of $H$, we have $c(xu_{1})\ne c(u_{1}u_{2})$ and $c(yv_{1})\ne c(v_{1}v_{2})$. Thus, $xyv_{1}v_{2}u_{2}u_{1}x$ is a properly colored $C_{6}$.

{\bf Case 1.2} $|E(H[U])|=1.$

Assume $uv\in E(H[U])$ and $c(uv)$ is starred at $u$. Then we have $c(xu)\ne c(uv)$. Also, $c(vz)\ne c(uv)$. Take a vertex $w\in U\setminus\{u, v\}$. Since $d^{c}(w)\ge 2$, we have $|E_{H}(w,\{x, y ,z\})|\ge 2$. At least one of $\{x, z\}$, say $x$, such that $c(wx)$ is starred at $w$ and $c(wx)\ne c(wy)$. Also, we have $c(wx)\ne c(ux)$. Thus $wxuvzyw$ is a properly colored $C_{6}$.

{\bf Case 1.3} $E(H[U])=\0.$

For all $v\in U$, since $d^{c}(v)\ge 2$, we have $|E_{H}(v, \{x, y, z\})|\ge 2$.  Notice that $|U|\ge n-3\ge 5.$ If there are three vertices in $U$, say $u_{1}, u_{2}, u_{3}\in U$, such that they have a  common neighborhood $\{x, z\}$ in $H$, then at least one of $\{u_{1}x, u_{1}z\}$, say $u_{1}x$, such that $c(u_{1}y)\ne c(u_{1}x)$. Also, at most one edge of $\{u_{2}x, u_{2}z, u_{3}x, u_{3}z\}$ has the same color as $c(u_{2}u_{3})$. Thus, at least one of $\{xu_{1}yzu_{3}u_{2}x, xu_{1}yzu_{2}u_{3}x\}$ is a properly colored $C_{6}$.

Now we assume that there are at least two vertices in $U$, say $u_{1}, u_{2}$, such that they have a common neighborhood $\{x, y\}$ or $\{y,z\}$, say $\{x, y\}$ in $H$.
If there is a vertex $u_{3}\in U\setminus \{u_{1}, u_{2}\}$ such that $u_{3}y, u_{3}z \in E(H)$,
we have $c(zx)\notin \{c(xu_{1}), c(xu_{2}), c(zu_{3})\}$ and at most one edge of $\{u_{1}x, u_{1}y, u_{2}x, u_{2}y\}$ has the same color as $c(u_{1}u_{2})$. Thus, at least one of $xu_{1}u_{2}yu_{3}zx$ and $ xu_{2}u_{1}yu_{3}zx$ is a properly colored $C_{6}$.
If there is a vertex $u_{3}\in U\setminus \{u_{1}, u_{2}\}$ such that $u_{3}x, u_{3}z \in E(H)$,
at least one of $xu_{1}u_{2}yzu_{3}x$ and $xu_{2}u_{1}yzu_{3}x$ is a properly colored $C_{6}$.
We may assume that $U$ has a common neighborhood $\{x, y\}$ in $H$. Take four distinct vertices $u_{1}, u_{2}, u_{3}, u_{4}\in U$. At most one edge of $\{u_{1}x, u_{1}y, u_{2}x, u_{2}y\}$ has the same color as $c(u_{1}u_{2})$ and at most one edge of $\{u_{3}x, u_{3}y, u_{4}x, u_{4}y\}$ has the same color as $c(u_{3}u_{4})$. Thus there is a properly colored $C_{6}$ in $\{u_{1}u_{2}, u_{3}u_{4}, xu_{i}, yu_{i}: 1\le i\le 4\}$.

{\bf Case 2.} $\Delta(G)=1.$

Note that if $G$ has three independent edges, then we can find a properly colored $C_{6}$.  Recall that $|E(G)|\ge n-6\ge 2$.  We have $n=8$ and $|E(G)|=2$. Let $E(G)=\{xy, zw\}$ and $U=V(K_{8})\setminus \{x, y, z, w\}=\{u_{1}, u_{2}, u_{3}, u_{4}\}.$

{\bf Case 2.1} There is an edge $u_{i}u_{j}$ such that $c(u_{i}u_{j})$ is starred at $u_{i}$, say $c(u_{1}u_{2}) $ is starred at $u_{1}$.

If there is one vertex in $\{x, y, z, w\}$, say $x$, such that $c(u_{1}x)\ne c(u_{1}u_{2})$, then $u_{1}xyzwu_{2}u_{1}$ is a properly colored $C_{6}$. We  assume that $c(u_{1}x)=c(u_{1}y)=c(u_{1}z)=c(u_{1}w)=c(u_{1}u_{2}).$ Since $d^{c}(u_{1})\ge 2$, we can assume that $c(u_{1}u_{3})$ is starred at $u_{1}$ and $c(u_{1}u_{3})\ne c(u_{1}u_{2}).$  Thus $u_{1}xyzwu_{3}u_{1}$ is a properly colored $C_{6}$.

{\bf Case 2.2} For all edge $u_{i}u_{j}$, $c(u_{i}u_{j})$ is not starred at $u_{i}$ or $u_{j}$.

Since $d^{c}(u_{1})\ge 2$ and $d^{c}(u_{2})\ge 2$, we can find two distinct vertices $v_{1}, v_{2}\in \{ x, y, z, w\}$ such that $c(u_{1}v_{1})$ is starred at $u_{1}$ and $c(u_{2}v_{2})$ is starred at $u_{2}$. If $v_{1}=x$ and $v_{2}=y$, then $u_{1}xzwyu_{2}u_{1}$ is a properly colored $C_{6}$. If $v_{1}=x$ and $v_{2}=z$, then $u_{1}xywzu_{2}u_{1}$ is a properly colored $C_{6}$.
\end{proof}


\vskip.2cm

\n{\large\bf 5. $K_{4}^{-}$ and $K_{2,3}$}
\vskip.2cm

In this section, we will prove Theorems 5 and  6. First, we will determine the exact value of  $pr(K_{n}, K_{4}^{-})$.

\noindent{\bf Theorem 5.} For $n\ge 4$, $pr(K_{n}, K_{4}^{-})= \lfloor\frac{3(n-1)}{2}\rfloor. $
\begin{proof}
{\bf The lower bound:} Consider an edge-coloring of $K_{n}$ as follows.
Take a triangle $C_{3}=xyz$ of $K_{n}$ and a maximum matching $M=\{x_{1}y_{1}, x_{2}y_{2}, \cdots, x_{\lfloor\frac{n-3}{2}]}y_{\lfloor\frac{n-3}{2}\rfloor}\}$ of $K_{n}-\{x,y,z\}$. There is one vertex $w$ in $V(K_{n})\setminus (V(M)\cup \{x,y,z\})$ when $n$ is even. For $1\le i\le \lfloor\frac{n-3}{2}\rfloor$, color all the edges of $\{ux_{i} :u\in\{x, y, z, x_{j}, y_{j}, 1\le j\le i-1\}\}$ with color $c_{1i}$ and all the edges of $\{uy_{i} :u\in\{x, y, z, x_{j}, y_{j}, 1\le j\le i-1\}\}$ with color $c_{2i}$. If $n$ is even, color all edges of $\{uw: u\in V(K_{n}-w)\}$ with a new color. Finally,  assign distinct new colors to all edges of $C_{3}\cup M$. In such a coloring, there is no properly colored $K_{4}^{-}$, and the number of colors is $\lfloor\frac{3(n-1)}{2}\rfloor$.

{\bf The upper bound:} We will prove that for any $\lfloor\frac{3n-1}{2}\rfloor$ edge-coloring of $K_{n}$, there is a properly colored $K_{4}^{-}$ by induction on $n$. The base case $n=4$ is trivial.
 Consider a $\lfloor\frac{3n-1}{2}\rfloor$ edge-coloring of $K_{n}$. If there is a vertex $v$ such that $d^{c}(v)\le1$, then $|C(K_{n}-v)|\ge \lfloor\frac{3n-1}{2}]-1\ge\lfloor\frac{3(n-1)-1}{2}\rfloor$, and there is a properly colored $K_{4}^{-}$ in $K_{n}-v$ by induction. We may assume that $d^{c}(v)\ge 2$ for all $v\in V(K_{n})$. Let $G$ be a subgraph of $K_{n}$, such that $e\in E(G)$ if and only if the color $c(e)$ appears only once in $K_{n}$. Since $d^{c}(v)\ge 2$ for all $v\in V(K_{n})$, we have $|E(G)|\ge 2n-\lfloor\frac{3n-1}{2}\rfloor=\lceil\frac{n+1}{2}\rceil$ which implies there is a path $P_{3}=xyz$ in $G$. By the construction of $G$, if $e=uv\in E(G)$, the $c(e)$ is starred at $u$ and $v$. We consider the following two cases.

{\bf Case 1.} $xz\notin E(G)$.

In this case, $c(xz)$ is not starred at $x$ or $z$, say $x$. Since $d^{c}(x)\ge 2$, there is a vertex $w\not\in\{ x, y, z\}$ such that $c(xw)$ is starred at $x$. Then $c(xz), c(yw)\notin \{c(xy), c(yz), c(xw)\}$ and $\{xy, yz, xz, xw, yw\}$ is a properly colored $K_{4}^{-}$.

{\bf Case 2.} $xz\in E(G)$.

In this case, we can assume $c(ux)=c(uy)=c(uz)$ for all $u\in V(K_{n})\setminus\{x, y ,z\}$; otherwise we easily have a properly colored copy of $K_{4}^{-}$ in $K_{n}[x, y, z, u]$.
Thus we have $$|C(K_{n}-\{x,y\})|\ge \left\lfloor\frac{3n-1}{2}\right\rfloor-3\ge\left\lfloor\frac{3(n-2)-1}{2}\right\rfloor$$ and there is a properly colored $K_{4}^{-}$ in $K_{n}-\{x, y\}$ by the induction hypothesis.
\end{proof}

Now we prove the lower bound and upper bound of $pr(K_{n}, K_{2,3})$. We conjecture that the exact value is closer to the lower bound.

\noindent{\bf Theorem 6.} For $n\ge 5$,  $\frac{7}{4}n+O(1)\le pr(K_{n}, K_{2,3})\le 2n-1.$

\begin{proof}
{\bf Lower bound:} Let $n=4k+r$, where $1\le r\le 4$. Set $V(K_n)=V_1\cup\ldots\cup V_k\cup V_{k+1}$ such that $V_i\cap V_j=\emptyset$ for $i\not=j$, $|V_i|=4$ for $1\le i\le k$ and $|V_{k+1}|=r$. We color the edges with endpoints in the same set with $6k+\binom{r}{2}$ distinct colors and color the remaining edges with $k$ addition colors $c_{1}, c_{2}, \ldots, c_{k}$ such that all edges with endpionts in $V_{i}$ and $V_{j}$ are colored with $c_{\min\{i, j\}}$, where $i\ne j$. The total colors are $\frac{7}{4}n+O(1)$ and there is no properly colored $K_{2,3}$.

{\bf The upper bound:} We will prove that for any $2n$ edge-coloring of $K_{n}$, there is a properly colored $K_{2,3}$ by induction on $n$. The base case $n=5$  are trivial. Consider a $2n$ edge-coloring of $K_{n}$. If there is a vertex $v$ such that $d^{c}(v)\le 2$, then $|C(K_{n}-v)|\ge 2n-2$ and there is a properly colored $K_{2,3}$ in $K_{n}-v$ by induction. We may assume that $d^{c}(v)\ge3$ for all $v\in V(K_{n})$. Let $G$ be a subgraph of $K_{n}$ where $e\in E(G)$ if and only if the color $c(e)$ appears only once in $K_{n}$. Since $d^{c}(v)\ge3$ for all $v\in V(K_{n})$, we have $|E(G)|\ge 3n-2n=n.$
Note that for $n\ge4, ex(n, P_{4})\le n$ where equality holds for the graph of disjoint copies of $C_{3}$ (see \cite{erdos0}). So we will consider the following two cases.

{\bf Case 1.} $G$  contains a $P_{4}=xyzw$.

If $G[V(P_{4})]\cong K_{4}$, then we can assume $c(ux)=c(uy)=c(uz)=c(uw)$ for all $u\in V(K_{n})\setminus\{x,y,z,w\};$ otherwise we easily have a properly colored copy of $K_{2,3}$. Therefore
$$|C(K_{n}-\{x, y, z\})|\ge 2n-6=2(n-3)$$
and there is a properly colored copy of $K_{2,3}$ in $K_{n}-\{x, y, z\}$ by the induction hypothesis.

Now we consider the case $G[V(P_{4})]\not\cong K_{4}$.
Since $d^{c}(x)\ge 3$ and $d^{c}(w)\ge 3$, there is a vertex $u\in V(K_{n})\setminus\{x,y,z,w\}$ such that $c(xu)$ or $c(wu)$, say $c(xu)$ is starred at $x$ and $c(xu)\notin \{c(xy), c(xw)\}$. Therefore, $\{xy, yz, zw, xw, xu, zu\}$ is a properly colored $K_{2,3}$.

{\bf Case 2.} $G$ is the graph of disjoint copies of $C_{3}$.

Take a triangle $T_{1}=xyzx$ of $G$. Since $d^{c}(x)\ge 3$, there is a vertex $u\in V(K_{n})\setminus\{x,y,z\}$ such that $c(xu)$ is starred at $x$ and $c(xu)\notin \{c(xy), c(xz)\}$. Suppose $u$ belong to the triangle $T_{2}=uvwu$ of $G$. Then $\{xy, xu, zy, zu, vy, vu\}$ is a properly colored $K_{2,3}$.
\end{proof}

~\\

\n{\large\bf 6. Open problems}
~\\

Although the topic of this paper has been proposed by Manoussakis, Spyratos, Tuza and Voigt \cite{Manoussakis} about twenty years ago, there are a few results about it. In this paper, we get the relationship of $pr(K_{n}, G)$ and $ex(n, \mathcal{G'})$ by Theorem 1.  Many problems on $pr(K_{n}, G)$ are worth being studied and we state four problems here.

\noindent{\bf Problem 1.}  Recall that Conjecture 2 is still open in the range $k\ge 7$.

\noindent{\bf Problem 2.} By Theorem 2, we have  $pr(K_{n}, K_{4})\ge ex(n, C_{4})+1$. Also, it is easy to see that $pr(K_{n}, K_{4})=O(n^{\frac{3}{2}})$.  It is natural to ask for the exact upper bound, i.e.
$$pr(K_{n}, K_{4})\le (1+o(1)) ex(n, C_{4}) ?$$

\noindent{\bf Problem 3.} By Theorem 6, we have $\frac{7}{4}n+O(1)\le pr(K_{n}, K_{2,3})\le 2n-1.$ It is natural to ask for the exact value of $pr(K_{n}, K_{2,3})$. Furthermore, it is interesting  to determine $pr(K_{n}, K_{s,t})$.

\noindent{\bf Problem 4.} Let $T_{k}$ be a tree of $k$ edges. The famous conjecture of Erd\H os and S\'os says that
$$ex(n, T_{k})\le \frac{(k-1)n}{2}.$$
Jiang and West \cite{Jiang2} conjectured that
$$ar(K_{n}, T_{k})\le \frac{(k-2)n}{2}+O(1).$$
Notice that $pr(K_{n}, T_{k})\le ar(K_{n}, T_{k})$ and $pr(K_{n}, K_{1,k})= ar(K_{n}, K_{1,k})=\frac{(k-2)n}{2}+O(1).$ It is natural to conjecture that
$$pr(K_{n}, T_{k})\le \frac{(k-2)n}{2}+O(1).$$
Furthermore, it is interesting to investigate the upper bound of $pr(K_{n}, T_{k})$ when $T_{k}\ne K_{1,k}$.


\end{document}